\documentclass[9pt,technote]{IEEEtran}
\usepackage[cmex10]{amsmath}
\usepackage{graphicx}
\usepackage{subfig}
\usepackage{amssymb}
\newtheorem{theorem}{Theorem}
\newtheorem{lemma}{Lemma}
\newtheorem{corollary}{Corollary}
\title{Nonlinear Predictor Feedback for Input-Affine Systems with Distributed Input Delays}
\author{Anton~Ponomarev\\
\thanks{Manuscript received January 15, 2015, revised June 6, 2015, accepted October 22, 2015.}
\thanks{A.~Ponomarev is with the Department of Control Theory, Saint Petersburg State University, Russia (e-mail: anton.pon.math@gmail.com).}
\thanks{}
\thanks{}
}

\begin{document}

\maketitle

\begin{abstract}
Prediction-based transformation is applied to control-affine systems with distributed input delays. Transformed system state is calculated as a prediction of the system's future response to the past input with future input set to zero. Stabilization of the new system leads to Lyapunov--Krasovskii proven stabilization of the original one. Conditions on the original system are: smooth linearly bounded open-loop vector field and smooth uniformly bounded input vectors. About the transformed system which turns out to be affine in the undelayed input but with input vectors dependent on the input history and system state, we assume existence of a linearly bounded stabilizing feedback and quadratically bounded control-Lyapunov function. If all assumptions hold globally, then achieved exponential stability is global, otherwise local. Analytical and numerical control design examples are provided.
\end{abstract}


\IEEEpubid{0000--0000/00\$00.00~\copyright~2015 IEEE}

\section{Notation}

The symbol $PC(T, X)$ stands for the space of piecewise continuous functions mapping $T\subset R$ into a Euclidean space $X$. The $L^2$ norm of $\varphi\in PC\big([-h,0), R^m\big)$ is $\Vert\varphi\Vert$, i.e.,
\begin{equation}
\Vert \varphi \Vert^2 =
	\int_{-h}^0 \Vert\varphi(\theta)\Vert^2 \,d\theta.
\end{equation}

Given $u\in PC\big([t-h,t), R^m\big)$, where $h>0$, let $u_t$ be a function defined as $u_t(\theta)=u(t+\theta)$ for all $\theta\in[-h,0)$.

$\mathcal{O}(R)$ is the closed $R$-ball about the origin in a normed space, specifically,
\begin{align}
x\in\mathcal{O}(R) &\Leftrightarrow \Vert x\Vert^2 \leq R^2, \\
(x,\varphi)\in\mathcal{O}(R) &\Leftrightarrow \Vert x\Vert^2 + \Vert\varphi\Vert^2 \leq R^2.
\end{align}

\section{Introduction}

\subsection{Problem statement}

Consider the system
\begin{multline} \label{plant}
\dot{x}(t) = f\big(x(t)\big) + B_0\big(x(t)\big) u(t) + B_1\big(x(t)\big) u(t-h) \\
	+ \int_{-h}^0 B_\text{int}\big(\theta, x(t)\big) u(t+\theta) \,d\theta,
\end{multline}
where $x\in R^n$, $u\in R^m$, $h>0$, and the following assumptions hold for some $M_f<\infty$ and $R\in(0,\infty]$:

1) regarding $f$, for all $x,x_0\in \mathcal{O}(R)$:
\begin{align}
\Vert f(x)\Vert &\leq M_f \Vert x\Vert, \label{bass2} \\
f(x) &= f(x_0) + A(x_0)(x-x_0) + o(x-x_0), \label{bass1} \\
A &\in C^0\big(\mathcal{O}(R), R^{n\times n}\big);
\end{align}

2) regarding $B_1$, for all $x,x_0\in \mathcal{O}(R)$:
\begin{align}
B_1 &\text{ is bounded on } \mathcal{O}(R), \\
B_1(x) &= B_1(x_0) + \Big( \mathcal{B}_1^1(x_0)(x-x_0), \dots, \mathcal{B}_1^m(x_0)(x-x_0) \Big) \notag \\
	&\phantom{=}\:+ o(x-x_0), \\
\mathcal{B}_1^i &\in C^0\big(\mathcal{O}(R), R^{n\times n}\big) \quad \forall i\in\overline{1,m};
\end{align}

3) regarding $B_\text{int}$, for all $x,x_0\in \mathcal{O}(R)$ and $\theta\in[-h,0]$:
\begin{align}
B_\text{int} &\text{ is bounded on } [-h,0]\times\mathcal{O}(R), \\
B_\text{int}(\cdot,x) &\in PC\big([-h,0], R^{n\times m}\big), \\
B_\text{int}(\theta,x) &= B_\text{int}(\theta,x_0) \notag \\
	&\phantom{=}\:+ \Big( \mathcal{B}_\text{int}^1(\theta,x_0)(x-x_0), \dots, \mathcal{B}_\text{int}^m(\theta,x_0)(x-x_0) \Big) \notag \\
	&\phantom{=}\:+ o(x-x_0), \\
\mathcal{B}_\text{int}^i(\theta,\cdot) &\in C^0\big(\mathcal{O}(R), R^{n\times n}\big) \quad \forall i\in\overline{1,m}.
\end{align}

The linear bound (\ref{bass2}) is assumed primarily to establish simple bounds on the system's solutions.

This paper is concerned with stabilization of the origin of (\ref{plant}), however, the approach presented here is applicable to time-variant systems with more point-wise and integral delays.

As for the practical importance of (\ref{plant}), such systems are used, e.g., for modeling of networked control systems where control delay value is unknown and varies rapidly \cite{goeb2010}.

There may be quite a few ways to stabilize (\ref{plant}). See, e.g., \cite{maze2013} where even more general time-variant systems with matrices $B_0$, $B_1$, and $B_\textrm{int}$ dependent on $x_t$ are tackled, the result being that sometimes a nonlinear feedback of $x(t)$ is stabilizing. However, we are concerned here with the so-called \textit{predictor feedback}. It got the name because, in a sense, it predicts the system's future behavior and how the control currently being chosen affects this behavior. Formally speaking, it is based on a state transformation which puts the system into a delay-free form. If the transformed system is stabilized with a feedback, it is possible to prove that the feedback will also stabilize the original system.

The predictor feedback theory is currently well-developed for linear systems with distributed input delays \cite{mani1979,kwon1980,arts1982} and nonlinear single-delay systems \cite{krst2010}. This paper extends the methodology to nonlinear systems with distributed input delays (\ref{plant}) for the first time.

\IEEEpubidadjcol

\subsection{Outline of the paper}

In Section \ref{S_preliminaries}, already known predictor feedbacks are summarized.

In Section \ref{S_reduction}, nonlinear predictor-based state transformation is defined for (\ref{plant}). The transformed system is derived in Theorem \ref{T_plant_reduced}.

In Section \ref{S_feedback}, an approach to stabilization of (\ref{plant}) is proposed (see Theorem \ref{T_main}). It is based on assumed existence of a quadratically bounded control-Lyapunov function and linearly bounded stabilizing feedback for the predictor-transformed system. Stability of the original system's closed loop is proved by Lyapunov--Krasovskii method.

In Section \ref{S_examples}, some illustrative problems are discussed.

\section{Predictor feedback overview} \label{S_preliminaries}

Let us recall a control methodology first designed for linear systems in \cite{mani1979,kwon1980,arts1982} and named there ``finite spectrum assignment'', ``receding horizon approach'', and ``model reduction''. We describe here a certain rehash of the same method based on the ``prediction'' concept and called ``predictor feedback''.

\subsection{Linear predictor feedback: one delay} \label{SS_intro_lin_single}

To begin with, consider the linear system with one discrete delay:
\begin{equation}
\dot{x}(t) = Ax(t) + B_1 u(t-h).
\end{equation}
The principles of predictor feedback in this case are:
\begin{itemize}
\item given the current state $x(t)$ and past input $u_t$, \textit{predictor} $Y\big(x(t), u_t\big)$ of the future state $x(t+h)$ is calculated using the variation of constants as
\begin{equation}
Y\big(x(t), u_t\big) = e^{Ah}x(t) + \int_{-h}^0 e^{-A\theta} B_1 u(t+\theta) \,d\theta;
\end{equation}
\item \textit{state transformation} is applied to the system:
\begin{equation}
y(t) = Y\big(x(t), u_t\big);
\end{equation}
\item \textit{predictor-transformed system} is delay-free:
\begin{equation}
\dot{y}(t) = Ay(t) + B_1 u(t);
\end{equation}
\item $y(t)$ is stabilized, if possible, with a feedback $u(t) = Ky(t)$;
\item the feedback for the original system is $u(t) = KY\big(x(t), u_t\big)$.
\end{itemize}

\subsection{Linear predictor feedback: distributed delay} \label{SS_intro_lin}

Suppose now that distributed delays are present:
\begin{equation} \label{plant_lin}
\dot{x}(t) = Ax(t) + B_0 u(t) + B_1 u(t-h) +
	\int_{-h}^0 B_\text{int}(\theta) u(t+\theta) \,d\theta.
\end{equation}
In this case \textit{literal} prediction of $x(t+h)$ at the current time $t$ would be challenging because it requires either \textit{knowledge} of the future input or \textit{assumptions} about the future input. The former leads to a non-causal and thus non-implementable feedback, the latter -- to an implicitly defined controller.

In \cite{mani1979,kwon1980,arts1982}, however, ``prediction with \textit{zero future input}'' is used. It is not the \textit{literal} prediction, i.e., actual future input is not required to be zero. The feedback is constructed as follows:
\begin{itemize}
\item \textit{predictor} $Y\big(x(t), u_t\big)$ is defined as a prediction of $x(t+h)$ with past input $u_t$ and future input set to zero, i.e., assuming $u(t+\theta) = 0$ for all $\theta\in[0,h]$; formally,
\begin{align}
Y\big(x(t), u_t\big) &= \xi(h), \\
\xi'(s) &= A\xi(s) + B_1 u(t+s-h) \nonumber\\
	&\phantom{=}\:+ \int_{-h}^{-s} B_\text{int}(\theta) u(t+s+\theta) \,d\theta,\\
\xi(0) &= x(t)
\end{align}
or, equivalently,
\begin{equation}
Y\big(x(t), u_t\big) = e^{Ah}x(t) + \int_{-h}^0 Q(\theta) u(t+\theta) \,d\theta,
\end{equation}
where
\begin{equation}
Q(\theta) = e^{-A\theta}B_1 + \int_{-h}^\theta e^{A(h-\theta+\tau)} B_\text{int}(\tau) \,d\tau;
\end{equation}
\item predictor-based \textit{state transformation} is applied:
\begin{equation} \label{pred_lin}
y(t)=Y\big(x(t), u_t\big);
\end{equation}
\item \textit{predictor-transformed system} is again delay-free:
\begin{equation} \label{pred_dyn_lin}
\dot{y}(t) = Ay(t) + \big( e^{Ah}B_0 + Q(0) \big) u(t);
\end{equation}
\item $y(t)$ is stabilized with $u(t) = Ky(t)$;
\item the final feedback is
\begin{equation} \label{feedback_lin}
u(t) = KY\big(x(t), u_t\big).
\end{equation}
\end{itemize}

Observe that the $y$-system (\ref{pred_dyn_lin}) is clearly exponentially stable under the feedback $u=Ky$ if suitable $K$ exists. Stability of (\ref{plant_lin}) under (\ref{feedback_lin}) is less obvious as the closed loop is a time-delay system. Lyapunov--Krasovskii analysis is highly desirable. The central problem in this analysis is to find a functional $v(x, \varphi)$ defined for all $x\in R^n$, $\varphi\in PC\big([-h,0), R^m\big)$ which is positive definite in terms of the norm $\Vert x\Vert^2 + \Vert\varphi\Vert^2$ and decreasing along the closed-loop solutions, i.e., when $x=x(t)$, $\varphi=u_t$. It is convenient to construct the functional starting with a Lyapunov function $y^T Vy$ for the closed-loop $y$-system. The functional may be then designed as
\begin{equation} \label{functional_lin}
v(x, \varphi) = Y^T(x,\varphi) V Y(x,\varphi)
	+ \int_{-h}^{\,0} e^{\sigma\theta} \Vert \varphi(\theta) \Vert^2 \,d\theta,
\end{equation}
where $\sigma>0$ (see \cite{pono2016} for proofs).

A lot of problems regarding linear predictor feedback have been addressed since its discovery including stability \cite{krst2008b,beki2011,maze2012}, robustness \cite{pono2016,mich2003,jank2009,kara2013}, delay-adaptive versions \cite{bres2010,bres2012}, and practical implementation issues \cite{palm1980,mond2002,mond2003,eves2003}. The approach is available for systems with state delays and an input delay \cite{khar2014a,khar2015}.

\subsection{Nonlinear predictor feedback} \label{SS_intro_nonlin}

Predictor feedback technique has been already expanded on nonlinear systems with one discrete delay in the input:
\begin{equation}
\dot{x}(t) = f\big(x(t), u(t-h)\big).
\end{equation}
The feedback was constructed and studied in \cite{krst2010} using the same idea as in the linear case: construct a predictor $Y\big(x(t), u_t\big)$ of $x(t+h)$; introduce new variable $y(t)=Y\big(x(t), u_t\big)$ which transforms the system to
\begin{equation}
\dot{y}(t) = f\big(y(t), u(t)\big)
\end{equation}
and stabilize it with a feedback $u = \kappa(y)$. An analysis is then conducted which shows that stability of the original system is achieved by $u(t) = \kappa\big(Y\big(x(t), u_t\big)\big)$.

As some examples of further development, let us mention the delay-adaptive variation of the feedback \cite{bres2014} and feedback for state-dependent delay \cite{beki2013}.

\section{Predictor-based state transformation}\label{S_reduction}

Following the linear theory described in Subsection \ref{SS_intro_lin}, we define predictor $Y\big(x(t), u_t\big)$ for (\ref{plant}) via prediction of $x(t+h)$ under \textit{zero future input}. The nonlinear version of (\ref{pred_lin}) is
\begin{equation} \label{reduct_nonlin}
y(t) = Y\big(x(t),u_t\big)
\end{equation}
with mapping $Y(x,\varphi)$ defined for $x\in R^n$, $\varphi\in PC\big([-h,0], R^m\big)$ by the system
\begin{align}
Y(x,\varphi) &= \xi(h), \label{reduct_Ydef}\\
\xi'(s) &= f\big(\xi(s)\big) + B_1\big(\xi(s)\big) \varphi(s-h) \notag \\
	&\phantom{=}\:+ \int_{-h}^{-s}B_\textrm{int}\big(\theta, \xi(s)\big)\varphi(s+\theta)\,d\theta, \label{reduct_ode}\\
\xi(0) &= x. \label{reduct_init}
\end{align}
The definition is algorithmic: to find $Y(x,\varphi)$, solve (\ref{reduct_ode}) with initial condition (\ref{reduct_init}) for $s\in[0,h]$.

There may be problems with this definition if $\xi(s)$ goes to infinity between $s=0$ and $s=h$. However, the following lemma shows that the mapping is defined at least in a zero neighborhood.

\begin{lemma} \label{L_Y}
$Y(x,\varphi)$ is defined globally if (\ref{bass2}) is satisfied for $R=\infty$ or, otherwise, locally for
\begin{equation} \label{neighborhoodx}
(x,\varphi) \in \mathcal{O}\left( \frac{R}{\rho} \right),
\end{equation}
where $\rho = \sqrt{2}e^{M_fh} \max\big\lbrace 1, h\max_{\theta\in[-h,0], x\in\mathcal{O}(R)} \Vert B_\textrm{int}(\theta, x) \Vert + \max_{x\in\mathcal{O}(R)}\Vert B_1(x)\Vert \big\rbrace$. Furthermore, for these $x$ and $\varphi$
\begin{gather}
\Vert x\Vert^2 \leq
	\rho^2\big( \Vert Y(x, \varphi)\Vert^2 + \Vert\varphi\Vert^2 \big), \label{xleqY} \\
\Vert Y(x, \varphi)\Vert^2 \leq
	\rho^2\big( \Vert x\Vert^2 + \Vert\varphi\Vert^2 \big). \label{Yleqx}
\end{gather}
\end{lemma}

\begin{IEEEproof}
Suppose that (\ref{bass2}) holds true along the solution of the problem (\ref{reduct_ode}), (\ref{reduct_init}). Applying the Gronwall--Bellman, Cauchy--Schwarz, and Young inequalities, one deduces
\begin{equation}
\Vert\xi(s)\Vert^2 \leq
	\rho^2\big( \Vert x\Vert^2 + \Vert\varphi\Vert^2 \big)
\end{equation}
for all $s\in[0,h]$. It leads to the conclusion that if (\ref{neighborhoodx}) is satisfied, then $\Vert\xi(s)\Vert\leq R$, so (\ref{bass2}) does hold along $\xi(s)$, and $\xi(s)$ does not go to infinity. Let $s=h$ to get (\ref{Yleqx}). To arrive at (\ref{xleqY}), use the same estimations when $s$ goes backwards from $h$ to $0$.
\end{IEEEproof}

Let us ignore the possibility of predictor nonexistence for a while (until Theorem \ref{T_main}) and derive the transformed system.

\begin{theorem} \label{T_plant_reduced}
System (\ref{plant}) after transformation (\ref{reduct_nonlin}) assumes the form
\begin{equation} \label{plant_reduced}
\dot{y}(t) = f\big(y(t)\big) + B\big(y(t), u_t\big)u(t)
\end{equation}
with $B(y,\varphi)$ defined for every $y$ and $\varphi$ implicitly by the equations
\begin{align}
B(y, \varphi) &= B_1\big(\xi(h)\big) + \beta(h), \label{B}\\
\beta'(s) &= \tilde{A}\big(s, \xi(s), \varphi\big) \beta(s) + B_\text{int}\big(-s, \xi(s)\big), \label{beta_ode}\\
\beta(0) &= B_0\big(\xi(0)\big), \label{beta_init}\\
\xi'(s) &= f\big(\xi(s)\big) + B_1\big(\xi(s)\big) \varphi(s-h) \notag \\
	&\phantom{=\;} +\int_{-h}^{-s} B_\textrm{int}\big(\theta, \xi(s)\big)\varphi(s+\theta)\,d\theta, \\
\xi(h) &= y, \\
\tilde{A}(s,\xi,\varphi) &= A(\xi) + \sum_{i=1}^m \mathcal{B}_1^i(\xi) \varphi_i(s-h) \notag \\
	&\phantom{=}\:+ \sum_{i=1}^m \int_{-h}^{-s} \mathcal{B}_\text{int}^i(\theta, \xi) \varphi_i(s+\theta) \,d\theta
\end{align}
granted that $\xi(s)$ exists on the interval $[0,h]$.
\end{theorem}

\begin{IEEEproof}
We are going to calculate the difference between $y(t)$ and $y(t+\Delta)$, where $\Delta$ is a small number, and find $\dot{y}(t)$ from there.

According to (\ref{reduct_nonlin}), $y(t)=\xi(h)$, where $\xi(s)$ satisfies
\begin{align}
\xi'(s) &= f\big(\xi(s)\big) + B_1\big(\xi(s)\big) u(t+s-h) \nonumber \\
	&\phantom{=}\:+ \int_{-h}^{-s}B_\text{int}\big(\theta, \xi(s)\big)u(t+s+\theta)\,d\theta, \label{xi_ode}\\
\xi(0) &= x(t). \label{xi_init}
\end{align}
Notice that this function $\xi$ is the same as in the statement of the theorem, just defined from the other end of the interval $[0,h]$. Likewise, $y(t+\Delta)=\tilde{\xi}(h)$, where $\tilde{\xi}(s)$ satisfies
\begin{align}
\tilde{\xi}'(s) &= f\big(\tilde{\xi}(s)\big) + B_1\big(\tilde{\xi}(s)\big) u(t+\Delta+s-h) \nonumber \\
	&\phantom{=}\:+ \int_{-h}^{-s}B_\text{int}\big(\theta, \tilde{\xi}(s)\big)u(t+\Delta+s+\theta)\,d\theta,\\
\tilde{\xi}(0) &= x(t+\Delta).
\end{align}

Consider now $\zeta(s)=\tilde{\xi}(s)-\xi(s+\Delta)$ for $s\in[0,h-\Delta]$. It follows from the definition of $\xi$ and $\tilde{\xi}$ that
\begin{equation}
\zeta(0) = x(t+\Delta)-\xi(\Delta) = B_0\big(\xi(0)\big)u(t)\Delta + o(\Delta), \label{zeta_init}
\end{equation}
and
\begin{align}
\zeta'(s) &= \tilde{\xi}'(s) - \xi'(s+\Delta) \nonumber \\
	&= f\big(\tilde{\xi}(s)\big) - f\big(\xi(s+\Delta)\big) \nonumber \\
	&\phantom{=}\:+ \Big(B_1\big(\tilde{\xi}(s)\big) - B_1\big(\xi(s+\Delta)\big)\Big)u(t+\Delta+s-h) \notag \\
	&\phantom{=}\:+ \int_{-h}^{-s} B_\text{int}\big(\theta,\tilde{\xi}(s)\big) u(t+\Delta+s+\theta)\,d\theta \notag \\
	&\phantom{=}\:- \int_{-h}^{-s-\Delta} B_\text{int}\big(\theta, \xi(s+\Delta)\big)u(t+\Delta+s+\theta)\,d\theta \nonumber \\
	&= A\big(\xi(s)\big)\zeta(s) + \sum_{i=1}^m \mathcal{B}_1^i\big(\xi(s)\big) \zeta(s) u_i(t+s-h) \notag\\
	&\phantom{=}\:+ \sum_{i=1}^m \int_{-h}^{-s} \mathcal{B}_\text{int}^i\big(\theta, \xi(s)\big) \zeta(s) u_i(t+s+\theta) \,d\theta \notag\\
	&\phantom{=}\:+ B_\text{int}\big(-s, \xi(s)\big)u(t)\Delta + o(\Delta) \notag\\
	&= \tilde{A}\big(s,\xi(s),u_t\big) \zeta(s) + B_\text{int}\big(-s, \xi(s)\big)u(t)\Delta + o(\Delta). \label{zeta_ode}
\end{align}
Equation (\ref{zeta_ode}) holds only if $\xi(s)\in\mathcal{O}(R)$. Observe that, up to $o(\Delta)$, the solution of the problem (\ref{zeta_ode}), (\ref{zeta_init}) coincides with that of the problem (\ref{beta_ode}), (\ref{beta_init}) multiplied by $u(t)\Delta$:
\begin{align}
\zeta(h-\Delta) &= \beta(h)u(t)\Delta + o(\Delta).
\end{align}
Therefore,
\begin{align}
\frac{y(t+\Delta)-y(t)}{\Delta} &= \frac{\tilde{\xi}(h)-\xi(h)}{\Delta} \nonumber \\
	&= \frac{\tilde{\xi}(h)-\tilde{\xi}(h-\Delta)+\zeta(h-\Delta)}{\Delta} \nonumber \\
	&= f\big(y(t)\big) + B\big(y(t),u_t\big)u(t) + \frac{o(\Delta)}{\Delta},
\end{align}
which proves the theorem.
\end{IEEEproof}

\section{Feedback design} \label{S_feedback}

The following result is inspired by the control-Lyapunov function concept.

\begin{theorem}\label{T_main}
If for all $(y,\varphi)\in\mathcal{O}(R)$ it is possible to define smooth functions $v_0(y)$ and $w_0(y)$ and a mapping $\kappa(y,\varphi)$ which satisfy the inequalities
\begin{gather}
m_{v_0}\Vert y\Vert^2 \leq v_0(y) \leq M_{v_0}\Vert y\Vert^2, \label{ineq_v0}\\
w_0(y) \geq m_{w_0}\Vert y\Vert^2, \label{ineq_w0}\\
\Vert \kappa(y,\varphi) \Vert \leq M_\kappa \Vert y\Vert, \label{ineq_kappa}\\
\big( f(y)+B(y,\varphi)\kappa(y,\varphi) \big)^T \nabla v_0(y) \leq
	-w_0(y), \label{lyapunoveq}
\end{gather}
where $m_{v_0}>0$ and $m_{w_0}>0$, then:
\begin{enumerate}
\item the origin is exponentially stable in the loop of (\ref{plant}) closed by the feedback
\begin{equation} \label{feedback}
u(t) = \kappa\big(Y(x, u_t), u_t\big);
\end{equation}
\item the stability is global if $R=\infty$ or local if $R<\infty$ with the region of attraction containing at least the ball
\begin{equation} \label{attraction}
\big(x(t), u_t\big) \in\mathcal{O}\left( \sqrt{\frac{m_v}{M_v}} \frac{R}{\rho} \right),
\end{equation}
where $\rho$ is given by Lemma \ref{L_Y}, and
\begin{align}
m_v &= \frac{1}{2\rho^2}\min\big\lbrace m_{v_0}, \gamma e^{-\sigma h},
	\rho^2 \gamma e^{-\sigma h} \big\rbrace,\\
M_v &= \gamma+\rho^2 M_{v_0},\\
\gamma &= \frac{m_{w_0}}{2M_\kappa^2},\\
\sigma &= \frac{m_{w_0}}{2M_{v_0}};
\end{align}
\item solutions of the closed loop (\ref{plant}), (\ref{feedback}) exhibit the exponential decay property
\begin{equation} \label{decay}
\Vert x(t)\Vert^2 + \Vert u_t\Vert^2 \leq
	\frac{M_v}{m_v} e^{-\sigma t} \big(
		\Vert x(0)\Vert^2 + \Vert u_0\Vert^2
	\big).
\end{equation}
\end{enumerate}
\end{theorem}

\begin{IEEEproof}
The proof is by the Lyapunov--Krasovskii method. Consider the functional for the closed loop (\ref{plant}), (\ref{feedback}):
\begin{equation} \label{functional}
v(x,\varphi) = v_0\big(Y(x,\varphi)\big) +
	\gamma \int_{-h}^0 e^{\sigma\theta} \Vert\varphi(\theta)\Vert^2 \,d\theta.
\end{equation}
Using (\ref{xleqY}), (\ref{Yleqx}), and (\ref{ineq_v0}), one can estimate
\begin{equation} \label{vestim}
m_v \big( \Vert x\Vert^2 + \Vert\varphi\Vert^2 \big) \leq
	v(x,\varphi)
\leq
	M_v \big( \Vert x\Vert^2 + \Vert\varphi\Vert^2 \big)
\end{equation}
with $m_v$ and $M_v$ given in the statement of the theorem.

Let $(x(t), u_t)$ be an arbitrary solution of the closed loop and suppose that (\ref{reduct_nonlin}) is defined along the whole trajectory $(x(t), u_t)$. Denote the value of $v(x,\varphi)$ along this solution as $v(t)$, i.e.,
\begin{equation}
v(t) = v\big(x(t), u_t\big).
\end{equation}
It follows from the choice of $\gamma$ and $\sigma$ and from (\ref{ineq_w0}), (\ref{ineq_kappa}), (\ref{lyapunoveq}) that
\begin{align}
\dot{v}(t) &\leq
	\big(-m_{w_0} + \gamma M_\kappa^2\big) \big\Vert Y\big(x(t), u_t\big)\big\Vert^2 \nonumber \\
&\phantom{\leq}\:- \sigma \gamma \int_{-h}^0 e^{\sigma\theta}
	\Vert u(t+\theta)\Vert^2 \,d\theta \nonumber \\
&\leq -\sigma v(t). \label{dv}
\end{align}

Combining (\ref{vestim}) and (\ref{dv}) leads to (\ref{decay}). Finally, (\ref{decay}) together with (\ref{attraction}) taken at $t=0$ ensures that at no time $t\geq 0$ does the system leave the region $\big(x(t), u_t\big)\in\mathcal{O}(R/\rho)$ which, according to Lemma~\ref{L_Y}, is sufficient for (\ref{reduct_nonlin}) to be defined for all $t\geq 0$. Thus, the proof just given is correct inside (\ref{attraction}).
\end{IEEEproof}

The next statement is supplemental to Theorem \ref{T_main}. It presents a way to construct the feeedback $\kappa(y,\varphi)$, given $v_0(y)$ and $w_0(y)$. Its application is demonstrated in Section \ref{SS_ex3}.

\begin{corollary} \label{C_default_feedback}
Define $\mathbb{B}\subset R^{n\times m}$ as the set of all possible values of $B(y,\varphi)$ when $(y,\varphi)\in\mathcal{O}(R)$. Suppose $\mathbb{B}$ is bounded.

If there are smooth functions $v_0(y)$ and $w_0(y)$ meeting the conditions (\ref{ineq_v0}), (\ref{ineq_w0}), and
\begin{gather}
\Vert \nabla v_0(y) \Vert \leq M_{\nabla v_0} \Vert y \Vert, \\
f^T(y)\nabla v_0(y) - k\Vert B^T\nabla v_0(y)\Vert^2 \leq -w_0(y) \label{ineq_k}
\end{gather}
for some constants $M_{\nabla v_0}>0$, $k$, all $B\in\mathbb{B}$, and all $y\in\mathcal{O}(R)$, then the feedback
\begin{equation} \label{default_feedback}
\kappa(y, \varphi) = -k B^T(y, \varphi) \nabla v_0(y)
\end{equation}
satisfies (\ref{ineq_kappa}) and (\ref{lyapunoveq}). In other words, the implications of Theorem~\ref{T_main} hold true for these $v_0$, $w_0$, and $\kappa$.
\end{corollary}

\begin{IEEEproof} (\ref{ineq_kappa}) is satisfied with
\begin{equation}
M_\kappa = |k| \max_{B\in\mathbb{B}} \Vert B\Vert M_{\nabla v_0}.
\end{equation}
(\ref{lyapunoveq}) turns into (\ref{ineq_k}) after substitution of (\ref{default_feedback}).
\end{IEEEproof}

Assuming infinitely fast computations, the control algorithm suggested by Corollary \ref{C_default_feedback} is this:
\begin{enumerate}
\item At time $t$, given $x(t)$ and $u_t$, solve (\ref{xi_ode}), (\ref{xi_init}) for $\xi(s)$.
\item Assign $y(t) := \xi(h)$.
\item Solve (\ref{beta_ode}), (\ref{beta_init}) for $\beta(h)$.
\item Assign $B\big(y(t), u_t\big) := B_1 + \beta(h)$.
\item Apply $u(t) = -k B^T\big(y(t), u_t\big) \nabla v_0\big(y(t)\big)$ and repeat.
\end{enumerate}

\section{Examples} \label{S_examples}

\subsection{Scalar case}

The most trivial examples are equations of the form
\begin{equation} \label{ex1}
\dot{x}(t) = f\big(x(t)\big) + b_0 u(t) + b_1 u(t-h) + \int_{-h}^0 b_\text{int}(\theta) u(t+\theta) \,d\theta
\end{equation}
with same-sign input coefficients, e.g., $b_0>0$, $b_1>0$, and $b_\text{int}>0$. In this case $B(y,\varphi)$ is a positive scalar separated from zero, so a possible feedback is (\ref{feedback}) with
\begin{equation}
\kappa(y, \varphi) = \frac{-f(y)-y}{B(y, \varphi)}.
\end{equation}

When input coefficients are of different signs, $B(y, \varphi)$ may or may not be zero. If it is zero for some $y$ and $\varphi$, then the origin of (\ref{ex1}) may or may not be stabilizable. The exact case is unclear because, unlike in the linear systems, here $B(y, \varphi)$ is not constant.

\subsection{Explicit prediction} \label{SS_ex2}

As explained in \cite{krst2010}, for some nonlinear systems with one input delay the predictor can be obtained by quadrature. The same holds true for distributed input delays. Let us give an example.

The system
\begin{equation} \label{ex2}
\left\lbrace\begin{array}{l}
	\dot{x}_1(t) = x_2^2(t) + u(t-h), \\
	\dot{x}_2(t) = x_2(t) + u(t)
\end{array}\right.
\end{equation}
allows explicit predictor transformation
\begin{align}
y_1 &= x_1 + \frac{e^{2h}-1}{2}x_2^2 + \int_{-h}^0 u(t+\theta)\,d\theta, \\
y_2 &= e^h x_2
\end{align}
which results in the system
\begin{equation} \label{ex2_reduced}
\left\lbrace\begin{array}{l}
	\dot{y}_1 = y_2^2 + \big(1+ \big(e^h-e^{-h}\big) y_2\big)u, \\
	\dot{y}_2 = y_2 + e^h u.
\end{array}\right.
\end{equation}

The origin of (\ref{ex2_reduced}) is globally asymptotically stabilizable. Indeed, another state and input transformation
\begin{align}
z_1 &= y_1 - e^{-h}y_2 + \frac{e^{-2h}-1}{2}y_2^2, \\
z_2 &= e^{-h}y_2, \\
u &= -2z_2 + \tilde{u}
\end{align}
puts (\ref{ex2_reduced}) into the cascade form
\begin{equation} \label{ex2_final}
\left\lbrace\begin{array}{l}
	\dot{z}_1 = z_2^2 - z_2, \\
	\dot{z}_2 = -z_2 + \tilde{u}.
\end{array}\right.
\end{equation}
Consider the positive definite function (the design is from \cite{jank1996})
\begin{equation}
V(z) = \big( z_1 + \tfrac{1}{2}z_2(z_2-2) \big)^2 + z_2^2.
\end{equation}
Its time-derivative along the solutions of (\ref{ex2_final})
\begin{equation}
\dot{V} = -2z_2^2 + \frac{\partial V}{\partial z_2} \tilde{u}
\end{equation}
shows that the feedback
\begin{equation}
\tilde{u} = -\frac{\partial V}{\partial z_2}
\end{equation}
is globally asymptotically stabilizing.

On a side note, the overall transformation from (\ref{ex2}) to (\ref{ex2_final}) happens to be linear:
\begin{align}
z_1 &= x_1 - x_2 + \int_{-h}^0 u(t+\theta) \,d\theta, \\
z_2 &= x_2, \\
u &= -2x_2 + \tilde{u}.
\end{align}

\textbf{Remark:} Returning from $z$ back through $y$ to $x$ leads to a rather complex feedback for the original system (\ref{ex2}). Does it achieve global asymptotic stability (GAS)? \textit{On the one hand}, transformation $x\mapsto y$ is well defined globally and transformed system is closed-loop GAS. \textit{On the other}, assumption (\ref{bass2}) and conditions (\ref{ineq_v0})--(\ref{ineq_kappa}) hold only locally if anywhere. Therefore, we are bound to conclude that our proofs do not warrant GAS of the original systems. This observation is valid not only in this example. It represents an inherent limitation of our proofs of stability. See Conclusions for further remarks.

\subsection{Numerical prediction} \label{SS_ex3}

Consider the nonlinear inverted pendulum described by the system
\begin{equation} \label{ex3}
\left\lbrace\begin{array}{l}
	\dot{x}_1(t) = x_2(t), \\
	\dot{x}_2(t) = \sin x_1(t) + u(t) + u(t-h),\quad h=\pi/4.
\end{array}\right.
\end{equation}
The design presented below is also possible for all $h\in[0,\pi/2]$.

In the above notation,
\begin{equation}
f(x) = \begin{pmatrix}
	x_2 \\
	\sin x_1
\end{pmatrix}, \quad
B_0 = B_1 = \begin{pmatrix}
	0 \\
	1
\end{pmatrix}.
\end{equation}
Assumptions (\ref{bass2}) and (\ref{bass1}) hold true with
\begin{equation}
A(x) = \begin{pmatrix}
	0 & 1 \\
	\cos x_1 & 0
\end{pmatrix}, \quad M_f = 1, \quad R = \infty.
\end{equation}

We are going to use Corollary \ref{C_default_feedback} to design the feedback function $\kappa$ and then simulate the closed loop (\ref{ex3}), (\ref{feedback}).

\emph{Step 1}: estimate $\mathbb{B}$ in order to apply Corollary \ref{C_default_feedback}. The problem (\ref{beta_ode}), (\ref{beta_init}) is
\begin{align}
\beta'(s) &= A\big(\xi(s)\big)\beta(s), \\
\beta(0) &= \begin{pmatrix}
	0 \\ 1
\end{pmatrix}.
\end{align}
From the inclusions
\begin{align}
\frac{d}{ds}\big( \beta_1^2+\beta_2^2 \big) 
&\in \big[0, 2\big( \beta_1^2+\beta_2^2 \big)\big], \\
\frac{d}{ds}\left( \frac{\beta_1}{\beta_2} \right) 
&\in \left[
	1-\left( \frac{\beta_1}{\beta_2} \right)^2,
	1+\left( \frac{\beta_1}{\beta_2} \right)^2
\right]
\end{align}
valid for $s\in[0, \pi/2)$ we find
\begin{align}
\beta_1^2(h)+\beta_2^2(h) &\in \Big[ 1, e^{2h} \Big], \\
\frac{\beta_1(h)}{\beta_2(h)} &\in \big[\tanh h, \tan h\big],
\end{align}
which means that $\beta(h)$ resides in a sector of a circular ring centered at the origin. It implies a rough estimation of $\mathbb{B}$ as the sector moved by $B_1$.

\begin{figure}[!t]
\centering
	\includegraphics[width=8.5cm]{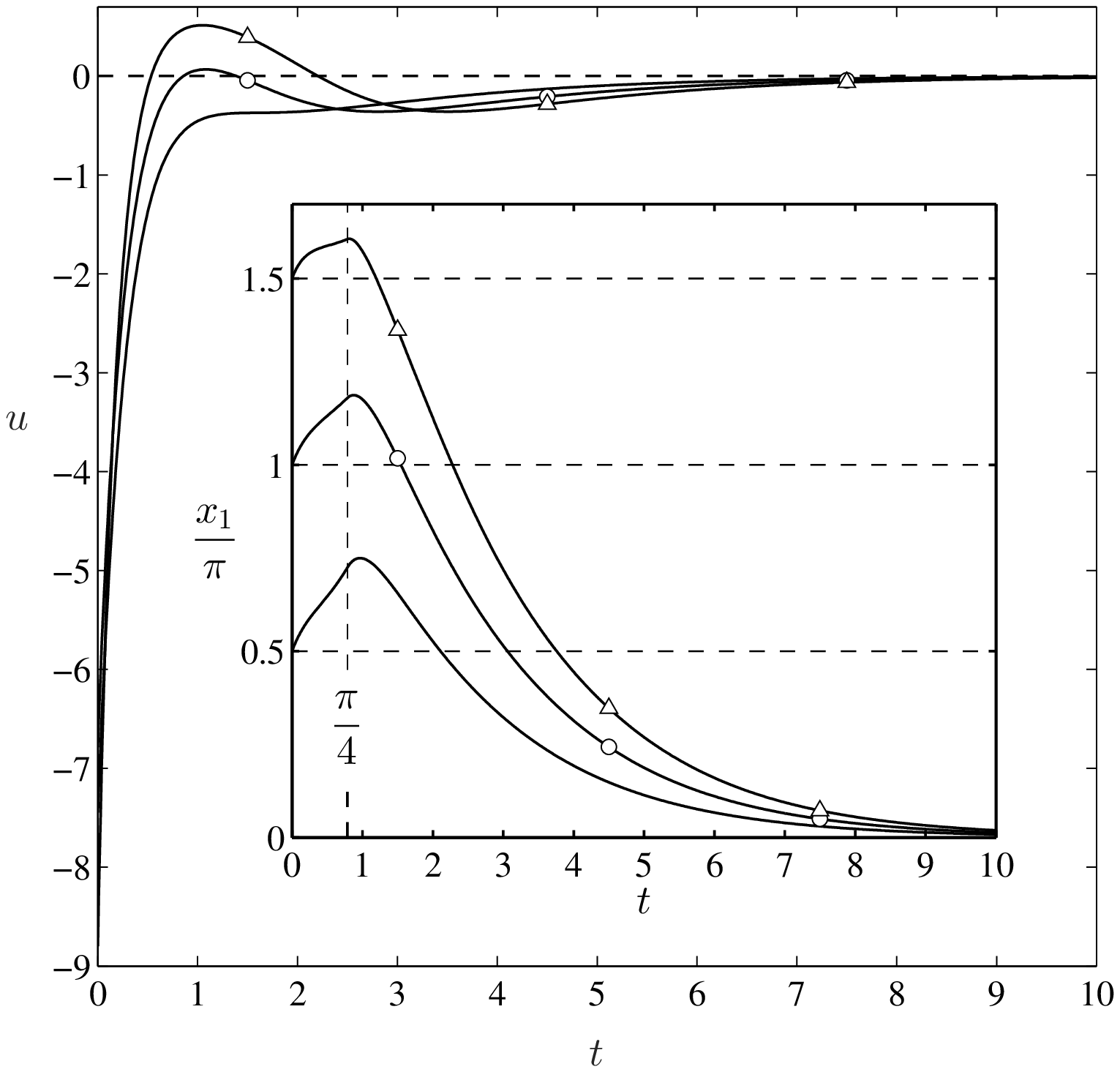}
\caption{Simulation of (\ref{ex3}) under the feedback (\ref{feedback}), (\ref{ex3_kappa}) with initial conditions $x_1(0)\in\lbrace \pi/2, \pi, 3\pi/2 \rbrace$, $x_2(0)=\pi/2$, $u_0\equiv 1$.}
\label{Figsim}
\end{figure}

\emph{Step 2}: find $v_0$, $w_0$ and $k$ that satisfy Corollary \ref{C_default_feedback}. Let
\begin{equation}
v_0(y) = y^T Vy
\end{equation}
with some positive definite $V$. The requirement of Corollary \ref{C_default_feedback} is then
\begin{equation} \label{ex3_lyap}
2f^T(y)Vy - 4k (Vy)^T BB^T Vy \leq -w_0(y).
\end{equation}
It is convenient to decompose $f(y)$ as
\begin{equation}
f(y) = F(y)y,
\end{equation}
where
\begin{gather}
F(y) = \begin{pmatrix}
	0 & 1 \\
	\alpha(y) & 0
\end{pmatrix}, \\
\alpha(y) = \frac{\sin y_1}{y_1} \in [-0.22, 1],
\end{gather}
then (\ref{ex3_lyap}) will follow from negative definiteness of the matrix
\begin{equation} \label{ex3_matrix}
V^{-1}F^T(y) + F(y)V^{-1} - 4k BB^T
\end{equation}
for all $y\in R^n$, all $B\in\mathbb{B}$, and some $k$. We found computationally that for $h=\pi/4$ it is valid to take $V=I$ and $k=1$, so the feedback function suggested by Corollary \ref{C_default_feedback} is
\begin{equation} \label{ex3_kappa}
\kappa(y, \varphi) = - B^T(y,\varphi) y.
\end{equation}

Fig. \ref{Figsim} shows the simulation results of (\ref{ex3}) under the feedback (\ref{feedback}), (\ref{ex3_kappa}) for different initial conditions. Euler's approximation with time step of $0.01$ was used for calculating $x(t)$, $\beta(s)$, and $\xi(s)$.

\section{Conclusions}

Let us pose some problems for a future discussion.

Firstly, concerning practical implementation of the proposed feedback, it should be mentioned that solving nonlinear equations may be costly, and any approximation, strictly speaking, requires robustness analysis.

Secondly, it is demonstrated in Section \ref{SS_ex2} that our linear-quadratic bounds (\ref{bass2}) and (\ref{ineq_v0})--(\ref{ineq_kappa}) may prohibit one from declaring global stabilization even when it is likely to be achieved. In fact, such restrictive bounds are not necessary. Their purpose is to simplify stability analysis by making the case ``almost linear''. The simplification is evinced in the exponential rate of decay (\ref{decay}) which is generally not expected from nonlinear systems. A possible way to avoid assumptions (\ref{bass2}) and (\ref{ineq_v0})--(\ref{ineq_kappa}) and to achieve class $\mathcal{KL}$ rate of decay in (\ref{decay}) is to follow the methodology of \cite{krst2010}.

\section*{Acknowledgment}

The author thanks Prof.~M.~Krstic for suggesting a number of substantial improvements to the presentation.

\bibliographystyle{IEEEtran}
\bibliography{IEEEabrv,nonlin_predictor_bib}

\end{document}